\documentclass{amsart}

\usepackage{amsmath, amsthm}
\usepackage[arc,all,knot,frame,curve,poly]{xy}
\usepackage{hyperref}

\numberwithin{equation}{section}

\newtheorem{Thm}{Theorem}[section]
\newtheorem*{Thma}{Theorem \ref{thm:hkr}}
\newtheorem*{Thmb}{Theorem \ref{thm:derham}}
\newtheorem{Prop}[Thm]{Proposition}
\newtheorem{Lem}[Thm]{Lemma}

\theoremstyle{remark}
\newtheorem{Rem}[Thm]{Remark}
\newtheorem*{Ack}{Acknowledgment}

\theoremstyle{definition}

\newcommand{\Z}{{\mathbb Z}}
\newcommand{\R}{{\mathbb R}}

\newcommand{\dd}{{\mathrm{d}}}

\newcommand{\calC}{\mathcal{C}}
\newcommand{\calD}{\mathcal{D}}

\newcommand{\calI}{\mathcal{I}}

\newcommand{\calL}{\mathcal{L}}

\newcommand{\calV}{\mathcal{V}}

\newcommand{\Hoch}{\mathsf{Hoch}}
\newcommand{\HHoch}{\mathsf{HHoch}}
\newcommand{\HKR}{\mathrm{HKR}}
\newcommand{\Coder}{\mathsf{Coder}}
\newcommand{\Hom}{\mathsf{Hom}}
\newcommand{\Der}{\mathsf{Der}}
\renewcommand{\div}{\mathrm{div}}

\begin{document}

\title[On the HKR map for graded manifolds]{On the 
Hochschild--Kostant--Rosenberg map for graded manifolds}

\author[A.~S.~Cattaneo]{Alberto S.~Cattaneo}
\address{Institut f\"ur Mathematik, Mathematisch-Naturwissenschaftliche 
Fakult\"at --- Universit\"at Z\"urich--Irchel ---
Winterthurerstrasse 190 --- CH-8057 Z\"urich --- Switzerland}
\email{alberto.cattaneo@math.unizh.ch}

\author[D.~Fiorenza]{Domenico~Fiorenza}
\address{Dipartimento di Matematica ``G. Castelnuovo'' --- 
Universit\`a degli Studi di Roma ``La Sapienza'' ---  Piazzale Aldo Moro, 2 
--- I-00185 Roma --- Italy}
\email{fiorenza@mat.uniroma1.it}

\author[R.~Longoni]{Riccardo~Longoni}
\address{Dipartimento di Matematica ``G. Castelnuovo'' --- 
Universit\`a degli Studi di Roma ``La Sapienza'' ---  Piazzale Aldo Moro, 2 
--- I-00185 Roma --- Italy}
\email{longoni@mat.uniroma1.it}

\thanks{A.~S.~C.\ acknowledges partial support of SNF Grant No.~20-100029/1, 
D.~F.\ acknowledges partial support of the European Commission through its 
6th Framework Programme ``Structuring the European Research Area'' and the 
contract Nr. RITA-CT-2004-505493 for the provision of Transnational Access 
implemented as Specific Support Action.}


\begin{abstract}
We show that the Hochschild--Kostant--Rosenberg map from the space of 
multivector fields on a graded manifold $N$ 
(endowed with a Berezinian volume) to the cohomology of the algebra
of multidifferential operators on $N$ (as a subalgebra of the Hochschild 
complex of $C^\infty(N)$) is an isomorphism of
Batalin--Vilkovisky algebras.  These results generalize to differential
graded manifolds. 
\end{abstract}

\maketitle

\section{Introduction}
\label{sec:intro}

The multivector fields on a smooth manifold $M$ can be seen as 
multidifferential operators on the algebra $\calC^\infty(M)$ of smooth 
functions on $M$. This assignment is a particular case of the following 
general construction: given a graded associative and commutative algebra
$A$, one defines  the Hochschild--Kostant--Rosenberg map
$$\HKR\colon\calV^\bullet (A)\to 
\Hoch^\bullet(A)$$ from the space of multivector fields $\calV^\bullet(A):=
S^{\bullet}(\Der(A)[-1])[1]$ to the Hochschild complex
$\Hoch^\bullet(A)$, as the map which regards a multiderivation of $A$
as a multilinear operator. Actually the image of HKR is contained
in the subcomplex $\calD^\bullet(A)\subset \Hoch^\bullet(A)$ of 
multidifferential operators. 

If one considers $\calV^\bullet(A)$ as a complex 
with trivial differential, then the HKR map is a morphism of complexes, and
the classical Hochschild--Kostant--Rosenberg Theorem \cite{HKR} states that 
when $A$ is a smooth algebra, e.g., a polynomial algebra, the HKR map 
induces isomorphisms in cohomology $\calV^\bullet(A)\simeq 
{\mathsf H}^\bullet(\calD^\bullet(A)) \simeq \HHoch^\bullet(A)$.
In this paper we are primarily concerned with the case in which $A$ is
the  algebra of smooth functions on a graded manifold $N$. In this case
it is  known that HKR still induces an isomorphism  
$\calV^\bullet(N)\simeq {\mathsf H}^\bullet
(\calD^\bullet(N))$, where we used the short-hand notations $\calV^\bullet(N)$
for $\calV^\bullet(\calC^\infty(N))$ and $\calD^\bullet(N)$ for 
$\calD^\bullet(\calC^\infty(N))$; for a proof, see \cite{V} in case $N$ is 
an ordinary manifold and \cite{CF} for the general case.

Many interesting algebraic structures can be defined on the objects
introduced  above. It is well known that $\calV^\bullet(A)$ and
$\HHoch^\bullet(A)$ are  Gerstenhaber algebras \cite{Gerst}, that
${\mathsf H}^\bullet (\calD^\bullet(A))$ is a sub-Gerstenhaber algebra
of $\HHoch^\bullet(A)$,  and that HKR preserves these structures.
Moreover, when $A$ is a finite dimensional algebra endowed with a
non-degenerate symmetric  inner product compatible with the
multiplication of $A$, then 
$\calV^\bullet(A)$, ${\mathsf H}^\bullet(\calD^\bullet(A))$ and 
$\HHoch^\bullet(A)$ become Batalin--Vilkovisky (BV) algebras \cite{Tr}.
The purpose of this paper is to extend this construction to the 
case in which $A$ is the algebra of smooth functions on a graded manifold $N$.
In this case the algebra is not finite dimensional but we can remedy when 
$N$ has a Berezinian volume. We prove in fact the following
\begin{Thma}
Let $N$ be a graded manifold endowed with a fixed Berezinian volume $v$ 
and whose body is a closed smooth manifold. Then $\calV^\bullet
(N)$ and ${\mathsf H}^\bullet(\calD^\bullet(N))$ 
can be endowed with BV algebra structures compatible with their classical 
Gerstenhaber structures. Moreover $\HKR$ is a map of BV algebras.
\end{Thma}
The BV algebra structure on multidifferential operators is inspired by 
\cite{Tr}, whereas the BV structure on $\calV^\bullet(N)$ is the
standard  one on the space of multivector fields of a graded manifold
$N$. Both  structures depend on the choice of a Berezinian volume on
$N$ \cite{KSM}.  The HKR map lifts to an $L_\infty$ map \cite{K,CF}
and, at least in the non graded case, to a $G_\infty$ map \cite{Ta}
between complexes.
One may conjecture that it also lifts to a $BV_\infty$ map \cite{TT}. 
This would be the analogue, for a graded manifold, of Kontsevich's
cyclic formality conjecture \cite{S}.\\

In the second part of the paper, we generalize our results to
differential graded manifolds $(N,Q)$. From an
algebraic point of view, this corresponds to considering differential
graded commutative associative algebras $(A,\dd)$. In this case, the
Hochschild complex is actually a bicomplex 
with differentials $\delta_0$ and $\delta_1$, 
and the Hochschild cohomology will be the cohomology of the total 
complex. The Hochschild bicomplex and its cohomology will be denoted by
$\Hoch_{\mathsf{DG}}^\bullet(A)$ and $\HHoch_{\mathsf{DG}}^\bullet(A)$
to distinguish them from the Hochschild complex and cohomology of
$A$ seen as a graded algebra. The differential $\dd$ gives rise 
to the differential $\{\dd,\cdot\}$ on the space
${\calV}^\bullet(A)$ of multivector fields; the HKR map
$({\calV}^\bullet(A),\{\dd,\cdot\},0)\to (\Hoch_{\mathsf{DG}}^\bullet(A),\delta_0,\delta_1)$
(see Lemma~\ref{lem:dr}) 
is a map of bicomplexes. We show by an example that the induced map
in cohomology is not an isomorphism in general. In particular  we
consider the differential graded manifold $N=T[1]M$, where
$M$ is a smooth manifold, with $\dd$ given by the de~Rham differential,  
so that
$\calC^\infty(T[1]M)$ is the de~Rham  algebra $\Omega^\bullet(M)$ of
$M$, and we prove the following
\begin{Thmb}
If $M$ is a simply connected closed oriented smooth manifold of
positive dimension, then the HKR map ${\mathsf H}^\bullet(\calV^\bullet
(\Omega^\bullet(M)),\{\dd,\cdot\})\to\HHoch_{\mathsf{DG}}^\bullet
(\Omega^\bullet(M))$ is not an isomorphism.
\end{Thmb}
The key ingredient of the proof is the isomorphism \cite{Chen} between the 
(shifted) homology $H_\bullet(\calL M)[\dim M]$ of the free loop space 
$\calL M$ of $M$ and the Hochschild cohomology of the differential graded 
algebra $\Omega^\bullet(M)$ . 
We remark that when only ordinary smooth manifolds are considered, it is 
not known whether the space of multivector fields is quasi-isomorphic to 
the Hochschild cohomology. Up to our knowledge, only a partial result in
this direction is known 
\cite{N}, namely, when $M$ is a smooth manifold,
$\calV^\bullet(M)$ is quasi-isomorphic to
the topological Hochschild complex $\HHoch_{\mathsf
{top}}(\calC^\infty(M))$  consisting of continuous multilinear
homomorphisms (with respect to the  Fr\'echet topology).

If we further assume that $(N,Q)$ is an SQ-manifold, i.e., that the
vector field
$Q$ is divergence-free, then a BV structure is induced on the cohomology
${\mathsf H}^\bullet(\calV^\bullet (A),\{\dd, \cdot\})$ and on the
Hochschild cohomology
$\HHoch_{\mathsf{DG}}^\bullet (A)$, and the HKR map is a morphism of BV
algebras (although, as remarked above, not an isomorphism in general).
An example is the de Rham algebra $(\Omega^\bullet(M),\dd)$ of a closed
manifold $M$. In this case, the BV structure on
$\HHoch_{\mathsf{DG}}^\bullet (\Omega^\bullet(M))$ corresponds to the one 
found in \cite{CS} on the homology of the free loop space \cite{CJ,M}, whereas
the BV structure on ${\mathsf H}^\bullet
(\calV^\bullet (\Omega^\bullet(M),\{\dd,\cdot\}))$ is the trivial one.

The plan of the paper is as follows. We begin by 
constructing the BV  structure on the space of multivector fields in 
Section~\ref{sec:mvf}. Next we recall some facts on 
Hochschild cohomology in Section~\ref{sec:hoch}. 
Then we discuss BV structures on the space of multidifferential 
operators in Section~\ref{sec:mdo},
and in Section~\ref{sec:hkr} we define the HKR map, 
describe its main properties, and prove Theorem~\ref{thm:hkr}. Finally
in Sections~\ref{sec:dgm} and \ref{sec:sq} we present a generalization of 
these results to the case of differential graded manifolds and prove
Theorem~\ref{thm:derham}.
\\
\begin{Ack}
We thank Thomas Tradler and the Referee for useful comments on a first draft 
of the paper. R.~L.\ thanks the Universit\"at Z\"urich--Irchel and D.~F.\ 
thanks the IH\'ES for their hospitality. 
\end{Ack}

\section{BV structure on multivector fields}
\label{sec:mvf}

Let $A$ be a graded commutative and associative algebra and let $\Der(A)
= 
\oplus_{j\in\Z} \Der^j(A)$ be the graded Lie algebra of derivations of $A$,
namely $\Der^j(A)$ consists of linear maps $\phi\colon A\to A$ of degree
$j$ such that $\phi(ab) = \phi(a)b+ (-1)^{j\, |a|} a\phi(b)$ and the bracket 
is $\{\phi,\psi\}=\phi\circ\psi -(-1)^{|\phi||\psi|}\psi\circ\phi$. 

The space of multiderivations 
$\calV^\bullet(A):=S^\bullet(\Der(A)[-1])[1]$ 
can be endowed with a Gerstenhaber structure, with the wedge product and 
the bracket which is the extension of the graded commutator
$\{\cdot,\cdot\}$  on $\Der(A)$ to $\calV^\bullet(A)$ by the Leibnitz
rule. Since $A$ is graded, the space $\calV^\bullet(A)$ has a natural
double grading given by 
\[
\calV^{i,j}(A)= \{\phi\in
S^i(\Der(A)[1])[-1]\,|\,\deg(\phi)=j\}.
\]
We want to construct an operator $\Delta$ on 
$\calV^\bullet(A)$ which makes this Gerstenhaber algebra into a BV 
algebra. We will use as an auxiliary tool the complex ${\mathcal
I}^\bullet(A)$ of integral forms of $A$, closely following
\cite{D}; a different approach to the BV 
algebra structures on $\calV^\bullet(A)$ can be found in
\cite{KSM}. Denote by $\Omega^1(A)$ the space of $1$-forms of $A$, namely,
the space ${\rm Hom}(\calV^1(A),A)$, and assume that the Berezinian ${\rm
Ber}(\Omega^1(A))$ is free and generated by one element $v$.
To a \emph{divergence operator} $\div$, viz.\ an
even linear map
$\div\colon \Der(A)\to A$ satisfying 
\[
\div(fX)=f\div(X)+(-1)^{|f||X|}X(f),
\]
we associate a linear operator 
$L\colon \calV^1(A) 
\otimes_A {\rm Ber}(\Omega^1(A)) \to {\rm Ber}(\Omega^1(A))$ by the rule
\[
L(X\otimes v)=\div(X)\,v.
\]
Observe that
for every $f\in A$ and every
$X\in\Der(A)$, we have $L_X(fv)=X(f)\,v + (-1)^{|f||X|}f\,L_X(v)$ where we 
are using the  notation $L_X(v):=L(X\otimes v)$. 

We now introduce the space ${\mathcal I}^\bullet(A)$ of integral forms 
\cite{D} as the
$A$-module generated by the elements of ${\rm
Ber}(\Omega^1(A))$ and by the operations $\iota_X$ with
$X\in\calV^1(A)$, acting on the left and subject to the rules
$[\iota_X,\iota_Y]=0$
and $\iota_{fX}=f\iota_X$. 
The action of $L_X$ is extended to ${\mathcal I}^\bullet(A)$ by the rule
$[L_X,\iota_Y]=\iota_{\{X,Y\}}$
One can define an exterior derivative $\dd$ on
${\mathcal I}^\bullet(A)$ by imposing
$\dd v=0$ and forcing Cartan's identity $\dd \iota_X+\iota_X\dd=L_X$.
Indeed, a consequence of Cartan's formula is that
$\dd(\iota_{X_1}\cdots\iota_{X_k}v)=L_{X_1}(\iota_{X_2}\cdots\iota_{X_k}
v)-\iota_{X_1}\dd(\iota_{X_2}\cdots\iota_{X_k}
v)$, and the action of $\dd$ on elements of ${\mathcal I}^\bullet(A)$ can
be computed inductively. The exterior derivative $\dd$ defined by this
procedure is a differential precisely when
$[L_X,L_Y]=L_{\{X,Y\}}$.  This is equivalent to the
vanishing of the curvature of $\div$; namely,
\[
\div(\{X,Y\}) - X(\div(Y)) + (-1)^{|X|\,|Y|}
Y(\div(X))=0.
\]
Once the generator $v$ of ${\rm
Ber}(\Omega^1(A))$ is fixed, iterated ``contractions'' $\iota_X$ induce an
isomorphism
\[
\calV^\bullet(A)\xrightarrow{\sim}{\mathcal I}^{\bullet}(A)
\]
and the differential $\dd$ induces on the space of multivector fields an
operator $\Delta$ of degree $-1$ such that $\Delta^2=0$. 
An easy computation shows that $\Delta(X)=\div(X)$ for any
 $X\in\Der(A)$, 
and that $\Delta$ satisfies the seven term relation
\begin{multline}
\label{eq:seven}
\Delta(a\wedge b\wedge c) + \Delta(a)\wedge b\wedge c + (-1)^{|a|} a\wedge
\Delta(b)\wedge c + (-1)^{|a|+|b|} a\wedge b\wedge \Delta(c) =\\ 
= \Delta(a\wedge b)\wedge c + (-1)^{|a|}a\wedge\Delta(b\wedge c) + 
(-1)^{(|a|+1)|b|}b\wedge \Delta(a\wedge c)
\end{multline}
and the compatibility with the bracket
\begin{equation}
\label{eq:defbv}
\{a,b\} := (-1)^{|a|}\left(\Delta(a\wedge b) - \Delta(a)\wedge b - 
(-1)^{|a|}a\wedge\Delta(b)\right).
\end{equation}
Therefore we have proved 
\begin{Lem}
If the Berezinian ${\rm Ber}(\Omega^1(A))$ is a free 
$A$-module of rank one
and $\div$ is a curvature-free divergence operator,
then the operator $\Delta$ defined as above endows $\calV^\bullet(A)$ 
with a BV structure compatible with the usual Gerstenhaber structure.
\end{Lem}
\par
The main example of this construction is when $A=\calC^\infty(N)$, $N$ 
being a graded manifold endowed with a Berezinian volume $v$. 
In this case the operators
$L_X$ and $\iota_X$ are just the classical Lie derivatives and contraction
operators, and the complex ${\mathcal I}^\bullet(N)$ is the
complex of integral forms of the graded manifold. 
Since the Berezinian is a line
bundle and $v$ is a nowhere  zero section, 
there exists an operator $\div$ defined uniquely by the equation $L_Y(v) = 
\div(Y)\,v$, which is indeed a divergence operator whose curvature 
vanishes. Observe that in the case when $N$ is an oriented smooth
manifold, this
amounts to choosing an ordinary volume form $v$. In the case when $N=T[1]M$,
with $M$ an oriented smooth manifold, there is a canonical Berezinian
volume $v$ characterized by
\[
\int_N \alpha\,v = \int_M \alpha,
\qquad\forall\alpha\in C^\infty(N) = \Omega^\bullet(M).
\]
\begin{Rem}\label{rem:Fourier}
The geometry of $T[1]M$ is closely related to the geometry of the formal
neighborhood of $M$ inside its cotangent bundle $T^*M$. Namely, the
Liouville volume form on $T^*M$ induces a curvature-free divergence
operator $\Delta$ on $\calV^\bullet(T^*M)$, which makes it a
BV algebra.  The algebra $A=\Gamma(S^\bullet TM)$ of smooth
functions on
$T^*M$ which are polynomial along  the fibers is a BV subalgebra of
$\calV(T^*M)$; it can be considered as the algebra of multivector
fields on $T^*M$ which are ``infinitesimal in the cotangent direction''.
As a consequence of the ``Fourier transform''
\cite{CF,R}, the Gerstenhaber algebras $\calV^\bullet(T[1]M)$ and
$\calV^\bullet(A)$ are isomorphic. But it can be easily verified that
they are also isomorphic as BV algebras.
\end{Rem}

\begin{Rem}
For a smooth manifold $M$, integral forms are
just ordinary differential forms and $\calI^\bullet(M)$ is naturally 
identified with $\Omega^\bullet(M)$.
On the other hand, for a graded manifold $N$ which is non trivial in odd 
degrees, the complex $\calI^\bullet(N)$ of integral forms is not isomorphic to
the de~Rham complex of $N$ (see \cite{D} for details).
\end{Rem}

\section{BV structure on Hochschild cohomology}
\label{sec:hoch}

The aim of this Section is to recall some standard facts about Hochschild 
cohomology and fix notations for the rest of the paper. We address the reader
to \cite{L} and \cite{Tr} for a comprehensive treatment.

\subsection{Hochschild cohomology}

Let $A=\oplus_{j\in\Z}A_j$ be a graded 
algebra over $\R$, with a graded commutative associative product $\mu$
and a unit
${\bf 1}$. 
We also suppose that $A$ is endowed with a non degenerate symmetric inner 
product compatible with the algebra multiplication, namely such that $\langle 
a,b\rangle = (-1)^{|a|\,|b|} \langle b,a\rangle$ and 
$\langle\mu(a\otimes b), c\rangle = \langle a, \mu(b\otimes c)
\rangle$. Finally, a graded bimodule $B$ over the algebra $A$ is given.

Let us set $T(A):=\bigoplus_{k\ge 0} A^{\otimes k}$ 
and $T^B(A) := \R\oplus \bigoplus_{k,l\ge 0} A^{\otimes k} 
\otimes B \otimes A^{\otimes l}$. 
It is well known that $T(A)$ is a coalgebra and $T^B(A)$ a bi-comodule over 
$T(A)$ with the coproducts
\begin{align*} 
T(A)& \to T(A)\otimes T(A)\\
(a_1,\ldots, a_n)&\mapsto \sum_{i=0}^{n}\left(a_1,\ldots,a_i\right) 
\otimes \left(a_{i+1},\ldots,a_n\right)
\end{align*}
and
\begin{align*}
T^B(A)&\to (T(A)\otimes T^B(A))\oplus (T^B(A)\otimes T(A)) \\
(a_1,\ldots,a_k,b,a_{k+1},\ldots,a_n)&\mapsto \sum_{i=0}^k(a_1,\ldots,a_i) 
\otimes (a_{i+1},\ldots,b,\ldots,a_n)+\\ &\phantom{mn} 
+ \sum_{i=k}^n(a_1,\ldots,b,\ldots,a_i) \otimes (a_{i+1},\ldots,a_n).
\end{align*}
Hence we can define the space $\Coder(T(A),T^B(A))$,
of coderivations from $T(A)$ to $T^B(A)$, with 
respects to the above coproducts.

The Hochschild cochain complex of $A$ with values in $B$ is defined as
\[
\Hoch^\bullet(A,B):= \Coder(T(A[1]), T^{B[1]}(A[1]))[-1]
\]
where by $A[1]$ we mean the graded algebra obtained by shifting the degrees of 
$A$ by 1; namely, $A[1] = \oplus_{j\in\Z}(A[1])_j$ with 
$(A[1])_j := A_{j+1}$.
As usual one can make the identification $$\Hoch^\bullet(A,B) = 
\prod_{n}\Hom(A[1]^{\otimes
n},B[1])[-1]=
\prod_{n}\Hom(A^{\otimes n},B)[-n].$$

Let us denote by $\widetilde{\mu^B}$ and $\widetilde{\mu}$  the lifts of
the bimodule structure $\mu^B \colon A\otimes B
\otimes A\to B$ and of the multiplication $\mu\colon A\otimes A\to A$
 to coderivations of $T(A[1])$
 with values in $T^{B[1]}(A[1])$.
Then, on the Hochschild cochain complex we can define a
degree 1 differential 
$\delta^B \colon \Hoch^\bullet(A,B) \to \Hoch^\bullet(A,B)$, by setting 
 $\delta^B(f):=\widetilde{\mu^B}\circ f -(-1)^{|f|} f\circ 
\widetilde{\mu}$. It is easy to check that $(\delta^B)^2=0$; the
cohomology of the Hochschild complex with respect to the differential 
$\delta^B$ is called Hochschild cohomology of $A$ with values in $B$ and
it is denoted by $\HHoch^\bullet(A,B)$. When $B=A$ with the canonical
bimodule structure we write $\HHoch^\bullet(A)$ for 
$\HHoch^\bullet(A,A)$; moreover $\delta^A$ is simply denoted by 
$\delta$.

\begin{Rem} Since $A$ and $B$ are graded objects, the Hochschild complex
$\Hoch(A,B)$ is a bigraded object: in the identification 
$\Hoch^\bullet(A,B) = 
\prod_{n}\Hom(A[1]^{\otimes
n},B[1])[-1]$, the horizontal
degree is provided by the number of $A$-factors, and
the vertical degree by the degree of the maps:
\[
\Hoch^{i,j}(A,B) = \{f\in\Hom(A[1]^{\otimes
i},B[1])[-1]|\, \deg(f)=j\}. 
\] 
The differential $\delta^B$ is a horizontal differential, since it 
increases the number of factors by one, leaving the degree of the 
maps unchanged. So one can think of the Hochschild complex as a 
bicomplex, with horizontal differential $\delta_1^B(f):=\widetilde{\mu^B}\circ f -(-1)^{|f|}
f\circ \widetilde{\mu}$ and trivial vertical differential
$\delta_0^B:=0$, and to consider $\delta^B$ as the total differential
$\delta^B=\delta_0^B+\delta_1^B$. We will come back to this point
of view when we will be discussing the Hochschild cohomology of
differential graded algebras in Section \ref{sec:sq}.
\end{Rem}

\subsection{Operations on the Hochschild cochain complex}

On the Hochschild co\-chain complex 
$\Hoch^\bullet(A)$ one can define various operations. 
First, there is a composition $f\circ g$ whose 
graded antisymmetrization $\{f,g\}:=f\circ g - (-1)^{|f|\, |g|} g\circ f$ 
gives rise to a graded odd Lie bracket of degree $+1$, 
also known as the Gerstenhaber 
bracket. Notice that the associativity of the product $\mu$ of $A$ 
is equivalent to 
$\{\widetilde\mu,\widetilde\mu\} = 0$, which immediately implies that
the  Hochschild differential $\delta(f) =
\{\widetilde\mu, f\}$  indeed squares to zero. Similar
relations holds for $\widetilde{\mu^B}$  and $\delta^B$.

Next, using the identification of $\Hoch^\bullet(A)$ with
$\prod_{n\ge 0} \Hom(A^{\otimes n}, A)[-n]$ we define a
product between $\phi\in \Hom(A^{\otimes k}, A)[-k]$ and 
$\psi\in \Hom(A^{\otimes l}, A)[-l]$ as 
$$(\phi\cup \psi) (a_1\otimes \cdots \otimes a_{k+l}):=
(-1)^\epsilon 
\mu(\phi(a_1\otimes \cdots \otimes a_k) \otimes \psi(a_{k+1}\otimes \cdots 
\otimes a_{k+l})),$$ 
where $\epsilon=l(|a_1|+\cdots+|a_k|+k)$.
This associative product is non-commuta\-tive but it gives rise to a
graded  commutative product in cohomology. The cup product and the
Gerstenhaber  bracket satisfy in cohomology the graded Leibnitz rule 
\[
\{a,b\cup c\} = \{a,b\}\cup c + (-1)^{(|a|+1)|b|}b\cup\{a,c\}.
\]
Therefore $(\HHoch^\bullet(A),\cup,\{\cdot,\cdot\})$ is
a Gerstenhaber algebra \cite{Gerst}. 
In addition, on the complex $\Hoch^\bullet(A,A^*)$ one has an 
operator $\beta$ given by the dual to Connes' $B$-operator
\cite{Co}. More explicitly, one defines $\beta 
\colon \Hoch^\bullet(A,A^*) \to \Hoch^{\bullet-1}(A,A^*)$ as
\[
(\beta(f)(a_1 ,\ldots,a_n ))(a_{n+1}) :=
\sum_{i=1}^{n+1}
(-1)^\epsilon
(f(a_i ,\ldots,a_{n+1},a_1,\ldots,a_{i-1}))({\bf 1})
\]
where ${\bf 1}$ is the unit of $A$ and $\epsilon = |f|+ |a_1|+\cdots + 
|a_{n+1}| + (|a_i |+\cdots + |a_n|)(|a_1|+\cdots +|a_{i-1}|)$. 

The inner product on $A$ gives rise to an injection $P\colon A\to A^*$ which 
is an $A$-bimodule map, and, by composing the Hochschild 
cochains with the injection $P$, one obtains an injective map $\wp\colon 
\Hoch^\bullet(A) \to \Hoch^\bullet(A,A^*)$. If moreover $\wp$ is a 
quasi-isomorphism, i.e., induces an isomorphism $H(\wp)$ in cohomology, 
then we can define an operator $\Delta_\beta$ of degree $-1$
on $\HHoch^\bullet(A)$ by setting 
$\Delta_\beta = H(\wp)^{-1}\circ\beta\circ H(\wp)$. 
As shown in \cite{Tr} (see also \cite{Men}), the operator $\Delta_\beta$ 
squares to zero in cohomology and is compatible with 
the Gerstenhaber structure on $\HHoch^\bullet(A)$ in the sense that 
(cf.\ equation~\eqref{eq:seven})
\begin{multline*}
\Delta_\beta(a\cup b\cup c) + \Delta_\beta(a)\cup b\cup c + (-1)^{|a|} a\cup 
\Delta_\beta(b)\cup c + (-1)^{|a|+|b|} a\cup b\cup \Delta_\beta(c) =\\ 
= \Delta_\beta(a\cup b)\cup c + (-1)^{|a|}a\cup\Delta_\beta(b\cup c) + 
(-1)^{(|a|+1)|b|}b\cup \Delta_\beta(a\cup c)
\end{multline*}
and (cf.\ equation~\eqref{eq:defbv})
\begin{equation*}
\{a,b\} = (-1)^{|a|}\left(\Delta_\beta(a\cup b) - \Delta_\beta(a)\cup b - 
(-1)^{|a|}a\cup\Delta_\beta(b)\right).
\end{equation*}
In other words $(\HHoch^\bullet(A),\cup,\{\cdot,\cdot\},
\Delta_\beta)$ is a BV algebra. Summing up, we have

\begin{Prop}
\label{prop:bv}
If the map $\wp\colon \Hoch^\bullet(A) \to 
\Hoch^\bullet(A,A^*)$ induced by the inner product of $A$ is a 
quasi-isomorphism, then 
$\HHoch^\bullet(A)$ is endowed with a BV algebra structure, compatible with 
its Gerstenhaber structure. 
\end{Prop}

A trivial example is when $A$ is finite dimensional, and hence $\wp$
is an isomorphism. A more interesting case is the algebra of functions
on a graded manifold $N$ endowed with a Berezinian volume $v$.
In this case the pairing is defined by
\begin{equation}\label{e:pairing}
\langle f_1,f_2\rangle=\int_N f_1 f_2\, v.
\end{equation}
In general, when $N$ is a graded manifold,
$\Hoch^\bullet(\calC^\infty(N))$ is  not necessarily quasi-isomorphic to
$\Hoch^\bullet(\calC^\infty(N), \calC^\infty(N)^*)$, and hence we do not
know whether we can define a BV structure on $\Hoch^\bullet(\calC^\infty(N))$. 
However we will see in Section~\ref{sec:mdo} that a version of 
Proposition~\ref{prop:bv} can be applied to a certain subcomplex of the 
Hochschild complex, namely to the subcomplex of multidifferential
operators.

\section{BV structure on multidifferential operators}
\label{sec:mdo}

The Hochschild complex of $A$ has a sub-Gerstenhaber algebra 
$\calD^\bullet(A)$ consisting of multidifferential operators, namely sums 
of cochains of the form $(a_1,\ldots,a_n)\mapsto \prod_{i=1}^n\phi_i(a_i)$ 
where $\phi_i$ are compositions of derivations. The bigrading on the
Hochschild complex induces a bigrading on the subalgebra of
multidifferential operators: \[
\calD^{i,j}(A):=\calD^\bullet(A)\cap
\Hoch^{i,j}(A).\]
We now want to discuss under 
which conditions the cohomology of ${\mathcal D}^{\bullet}(A)$ admits a 
natural BV structure. As above we are assuming that there 
exists a non degenerate symmetric inner product on $A$ compatible with the 
multiplication, and hence an injective 
map $\wp\colon \Hoch^\bullet(A) \to \Hoch^\bullet(A,A^*)$. The point is to 
determine when the Connes cyclic $\beta$-operator $\beta\colon 
\Hoch^\bullet(A,A^*) \to \Hoch^{\bullet-1} (A,A^*)$ induces an operator 
$\Delta_\beta\colon{\mathcal D}^{\bullet}(A)\to {\mathcal D}^{\bullet-1}(A)$
making the diagram
\[
\xymatrix{
   {\mathcal D}^\bullet(A) \ar[r]^\wp
\ar@{-->}[d]_{\Delta_\beta} &  
\Hoch^\bullet(A,A^*) \ar[d]_\beta\\
 {\mathcal D}^{\bullet-1}(A)  \ar[r]^\wp & 
\Hoch^{\bullet-1}(A,A^*)\\
 }
\]
commutative. To answer this question, we look at the problem from a more 
general perspective; namely, let 
$C^\bullet(A)$ be any sub-Gerstenhaber algebra of
$\Hoch^\bullet(A)$  whose $\wp$-image in $\Hoch^\bullet(A,A^*)$ is closed under
$\beta$. Since $\wp$ is injective, $\beta$ induces a well-defined
operator $\Delta_\beta$ on the complex $C^\bullet(A)$. 
Following \cite{Tr} and \cite{Men}, the operator $\Delta_\beta$ 
squares to zero 
in the cohomology of $C^\bullet(A)$, and endows 
${\mathsf H}^\bullet(C^\bullet(A))$ with a BV algebra 
structure compatible with its Gerstenhaber 
structure. 
\par
We now specialize to the case when $A=C^\infty(N)$, where $N$ is a 
graded manifold endowed with a Berezinian volume $v$. 
In order to prove that the 
cohomology ${\mathsf H}^\bullet({\mathcal
D}^\bullet(N))$ of the algebra of multidifferential operator admits a
natural BV structure, we only need to prove that 
$(\beta\circ\wp)({\mathcal D}^\bullet(N))\subseteq
\wp({\mathcal D}^\bullet(N))$
with $\wp$ induced by the pairing \eqref{e:pairing}.
We first need the following ``integration-by-parts'' Lemma.
\begin{Lem}\label{l:multi}
Let $D$ be a multidifferential operator. Then there exist a
multidifferential operator $\tilde{D}$ such that
\[
\langle D(f_1,\dots,f_n),{\bf 1}\rangle = \langle
\tilde{D}(f_1,\dots,f_{n-1}), f_n\rangle
\]
\end{Lem}
Then we observe that for every $D\in {\mathcal D}^n(N)$ and for every
$i=1,\dots,n$, the operator
\[
D_i(f_1,\dots,f_n):= D(f_i,\dots,f_n,f_1,\dots,f_{i-1}),
\qquad f_1,\dots,f_n\in A,
\]
is still in ${\mathcal D}^n(N)$.
Finally
\begin{multline*}
(\beta\circ\wp(D))
(f_1 ,\ldots,f_{n-1})(f_{n}) =
\sum_{i=1}^{n}
(-1)^\epsilon
\langle D(f_i ,\ldots,f_{n},f_1,\ldots,f_{i-1}),{\bf 1}\rangle=\\
=\sum_{i=1}^{n}
(-1)^\epsilon
\langle D_i(f_1 ,\ldots,f_{n}),{\bf 1}\rangle=
\sum_{i=1}^{n}
(-1)^\epsilon
\langle \tilde D_i(f_1 ,\ldots,f_{n-1}),{f_n}\rangle=\\
=\wp\left(\sum_{i=1}^{n}
(-1)^\epsilon\tilde D_i\right)(f_1,\dots,f_{n-1})(f_n).
\end{multline*}

\begin{proof}[Proof of Lemma \ref{l:multi}]
The proof is by induction on the order of the multidifferential 
operator $D$. If $D$ is homogeneous of order zero, 
\[
D(f_1,\dots,f_{n})=\lambda f_1\cdots f_{n}
\]
for some constant $\lambda$, so that
\[
\langle D(f_1,\dots,f_{n}), {\bf 1}\rangle =
\int_N \lambda f_1\cdots f_{n}\, v=
\langle \lambda f_1\cdots f_{n-1},f_n\rangle
\]
and we are done. Now assume the claim
proved for operators up to order $k$ and prove it for order $k+1$
operators by the following argument. A homogeneous component of an order
$k+1$ multidifferential operator can be written as 
\[
D(f_1,\dots,f_{n})=D_0(f_1,\dots f_{i-1},X(f_i),f_{i+1},
\dots,f_{n})
\]
for a suitable multidifferential operator $D_0$ of order $k$, some index
$i$ and some vector field $X$. 
 We compute
\[
\langle D(f_1,\dots,f_n),{\bf 1} \rangle=
\langle D_0(f_1,\dots,X(f_i), \dots, f_{n}),{\bf 1}\rangle
\]
Here we have to distinguish two cases. If $i\neq n$, by the induction
hypothesis applied to $D_0$, we can write
\[
\langle D_0(f_1,\dots,X(f_i), \dots, f_{n}),{\bf 1}\rangle
=
\langle\tilde{D}_0(f_1,\dots,X(f_i), \dots, f_{n-1}), f_n\rangle
\]
and we are done. If $i=n$ then the induction hypothesis gives
\[
\langle D_0(f_1,\dots, f_{n-1},X(f_{n})), {\bf 1}\rangle
=
\langle\tilde{D}_0(f_1,\dots,f_{n-1}), X(f_n)\rangle.
\]
For any vector field $Y$, Cartan's formula gives 
$L_Y(v) = \dd i_Y(v)  +  i_Y \dd(v)=\dd i_Y(v) $, 
since $\dd(v)=0$ \cite{D}. Hence, by Stokes' Theorem we have that
\[
0=\int_N \dd i_Y (f \, v) = \int_N Y(f)\, v + (-1)^{|f|\, |Y|}
\int_N f\, L_Y(v).
\]
Recall for Section \ref{sec:mvf} that  
there exists an operator $\div$ defined uniquely by the equation $L_Y(v) = 
\div(Y)\,v.$ Therefore
\begin{equation}
\label{eq:sc}
\langle Y(f),{\bf 1}\rangle =
\int_N Y(f)\, v = - (-1)^{|f|\, |Y|} \int_N f\, \div(Y)\, v = -
\langle \div(Y),f\rangle.
\end{equation}

Going back to our problem with $D_0$, we apply the previous formula to the
vector field $Y=\tilde{D}_0(f_1,\dots,f_{n-1})X$ and obtain
\begin{align*}
\langle \tilde{D}_0(f_1,\dots,f_{n-1}),X(f_n)\rangle
&=\int_N \tilde{D}_0(f_1,\dots,f_{n-1})X(f_n)\, v\\
&=\langle \div(\tilde{D}_0(f_1,\dots,f_{n-1})X),f_n \rangle.
\end{align*}
The map $(f_1,\cdots,f_{n-1})
\mapsto \div(\tilde{D}_0(f_1,\dots,f_{n-1})X)$ is a multidifferential 
operator, and the Lemma is proved by setting 
$\tilde D (f_1,\dots,f_{n-1})=\div(\tilde{D}_0(f_1,\dots,f_{n-1})X)$.
\end{proof}

\section{The Hochschild--Kostant--Rosenberg map}
\label{sec:hkr}
The 
Hochschild--Kostant--Rosenberg  (HKR) map is defined as follows:
\begin{equation}
\label{eq:hkr-m}
\begin{array}{ccc}
\calV^\bullet(A) & \longrightarrow & \Hoch^\bullet(A)\\
\phi_1\wedge\cdots\wedge\phi_n & \mapsto &
\displaystyle{\frac1{n!}\sum_{\sigma\in S_n} 
\mathrm{sign}(\sigma)\ \phi_{\sigma(1)}\cup\cdots\cup\phi_{\sigma(n)}}.
\end{array}
\end{equation}
Note that the HKR map is actually a map of bigraded vector spaces:
$\calV^{i,j}(A)\to \Hoch^{i,j}(A)$. We have already observed that both
$\calV^\bullet(A)$ and
$\HHoch^\bullet(A)$ are Gerstenhaber algebras, and it is well known
that the HKR map in fact  preserves these structures. More explicitly
\begin{Thm}
\label{thm:hkrg} 
If $\calV^\bullet(A)$ is endowed with the zero differential, then $\HKR$ is a 
morphism of complexes. Moreover the induced map in cohomology is a
morphism of Gerstenhaber algebras.
\end{Thm}

\begin{proof}
This is a standard result: the fact that 
$\HKR$ respects the product structures in cohomology follows directly 
from the fact that the cup product is commutative in cohomology \cite{Gerst}. 
An easy check shows that for $X,Y\in\Der(A)$ we have
\[
\{\HKR(X), \HKR(Y)\} - \HKR(\{X,Y\}) = 0 
\]
and hence, by the compatibility between the bracket and the product, 
$\HKR$ induces in cohomology a map of Gerstenhaber algebras.
\end{proof}

The classical Theorem of Hochschild, Kostant and Rosenberg \cite{HKR}
states that when $A$ is a smooth algebra (e.g. for the coordinate ring
of a smooth affine algebraic variety) then the HKR map is a
quasi-isomorphism, i.e., induces an isomorphism
$\calV^\bullet(A)\xrightarrow{\sim}\HHoch^\bullet(A)$. 
\par
One sees from equation (\ref{eq:hkr-m}) that the HKR map actually takes
its values in the subcomplex
$\calD^\bullet(A)$ of multidifferential operators. 
For a smooth algebra $A$, the inclusion
$\calD^\bullet(A)\hookrightarrow \Hoch^\bullet(A)$ is a
quasi-isomorphism, so the classical Hochschild-Kostant-Rosenberg
theorem can then be stated as follows.
\begin{Thm}
\label{thm:hkr-iso}
If $A$ is a smooth algebra, then 
$\HKR\colon\calV^\bullet(A)\to 
{\mathsf H}^\bullet(\calD^\bullet(A))$ is
an isomorphism of Gerstenhaber algebras.
\end{Thm}
Our main result is a version of Theorem \ref{thm:hkr-iso} for graded manifolds, 
namely, we prove
\begin{Thm}
\label{thm:hkr}
Let $N$ be a graded manifold endowed with a fixed Berezinian volume $v$ 
and whose body is a smooth closed manifold. Then $\calV^\bullet
(N)$ and ${\mathsf H}^\bullet(\calD^\bullet(N))$ 
can be endowed with BV algebra structures compatible with their classical 
Gerstenhaber structures. Moreover 
 $\HKR\colon\calV^\bullet(N)\to 
{\mathsf H}^\bullet(\calD^\bullet(N))$ is
an  isomorphism of BV algebras.
\end{Thm}

\begin{proof}

We have seen in Sections~\ref{sec:mvf} and~\ref{sec:mdo} that, in case 
$A=\calC^\infty(N)$ is the 
algebra of smooth functions of a graded manifold $N$ endowed with a Berezinian 
volume form, then both $\calV^\bullet(N)$ and ${\mathsf H}^\bullet(
\calD^\bullet(N))$ are BV 
algebras in a way compatible with their classical Gerstenhaber structures.

We know from Theorem~\ref{thm:hkrg} that $\HKR$ induces 
in cohomology a morphism of Gerstenhaber algebras. Moreover we know from 
\cite{CF} that $\HKR\colon \calV^\bullet(N) \to \calD^\bullet(N)$ is a 
quasi-isomorphism. Therefore, by the compatibility between the BV Laplacian and
the Gerstenhaber bracket, we only need 
to prove that for every vector field $X\in\calV^1(N)$ on a graded manifold $N$,
we have 
$$\HKR(\Delta(X)) = \Delta_\beta(\HKR(X)).$$
To see this, consider the diagram
\[
\xymatrix{
\calV^1(N) \ar[r]^{\HKR} \ar[d]_{\Delta} & 
\calD^1(N) \ar[r]^{\wp\phantom{mmmmm}} \ar[d]_{\Delta_\beta}&  
\Hoch^1(\calC^\infty(N),\calC^\infty(N)^*) \ar[d]_\beta\\
\calV^0(N) \ar[r]^{\HKR} & 
\calD^0(N) \ar[r]^{\wp\phantom{mmmmm}} & 
\Hoch^0(\calC^\infty(N),\calC^\infty(N)^*)
}
\]
Since the diagram on the right commutes and $\wp$ is injective, 
commutativity of the diagram on the 
left follows from the commutativity of the external diagram. This is indeed 
the case since on the one side, 
for $X\in\calV^1(N)$ and $f\in\calC^\infty(N)$, we have that
\begin{equation}
\label{eq1}
\left(\beta(\wp(\HKR(X)))\right)(f) = - \langle X(f), {\bf 1}\rangle, 
\end{equation}
on the other side 
\begin{equation}
\label{eq2}
\left(\wp(\HKR(\Delta(X)))\right)(f) = \langle \Delta(X), f\rangle. 
\end{equation}
By Section~\ref{sec:mvf}, $\Delta(X)=\div(X)$,
and the right-hand sides of equations~\eqref{eq1} and \eqref{eq2} coincide 
by means of equation~\ref{eq:sc}.
\end{proof}

\section{The HKR theorem for differential graded
manifolds}\label{sec:dgm}

We now consider the more general case of differential graded manifolds,
i.e., of graded manifolds $N$ endowed with a degree 1 integrable vector
field $Q$. Note that, since the degree of $Q$ is 1, the integrability
condition $\{Q,Q\}=0$ is equivalent to $Q^2=0$. The algebraic counterpart
of a differential graded manifold $(N,Q)$ is a differential graded
algebra $(A,\dd)$, where $\dd$ is a degree one differential on the
graded algebra $A$. A classical example is given by the
de Rham algebra
$(\Omega^\bullet(M),\dd)$ of a differential manifold $M$ with the de
Rham differential. The corresponding graded manifold is 
$T[1]M$; the de Rham differential on differential forms corresponds to
a degree 1 integrable vector field on $T[1]M$. Note that ordinary graded
manifolds can be considered as differential graded manifolds with
the trivial vector field $Q=0$.
\par
The construction of the Hochschild
complex of a graded algebra $A$ with values in $B$ described in Section
\ref{sec:hoch} generalizes to the case of a differential graded 
algebra $(A,\dd)$. In this case one actually gets a nontrivial vertical
differential by setting $\delta^B_0(f):=
\widetilde\dd\circ f -(-1)^{|f|} f\circ\widetilde\dd$, where
$\widetilde\dd$ denotes the lift of the differential $\dd\colon A\to
A$, to coderivations of $T(A[1])$ with values in $T^{B[1]}(A[1])$.
The horizontal differential $\delta_1^B$ is the same as in the case of
graded algebras described in Section \ref{sec:hoch}. One easily checks
that the total differential $\delta^B=\delta_0^B+\delta_1^B$
squares to zero. We show this in the particular case $B=A$, the
general case being similar. By definition,
$\delta_1=\{\widetilde\mu,\cdot\}$ and
$\delta_0=\{\widetilde\dd,\cdot\}$; the associativity of the
product $\mu$ is equivalent to 
$\{\widetilde\mu,\widetilde\mu\} = 0$, 
the fact that $\dd$ is a derivation for $\mu$ is equivalent to
$\{\widetilde\dd,\widetilde\mu\} = 0$,
and $\dd^2=0$ is equivalent to 
$\{\widetilde \dd,\widetilde \dd\}=0$.
These three properties immediately imply that the 
Hochschild differential $\delta(f) = \{\widetilde\mu+\widetilde\dd,
f\}$ indeed squares to zero. 

The total complex will be denoted by $\Hoch_{\mathsf{DG}}(A,B)$; its
cohomology is called Hochschild cohomology of $A$ with values in
$B$ and it is denoted by $\HHoch_{\mathsf {DG} }^\bullet(A,B)$, where the
subscript 
${\mathsf {DG} }$ means that we are working in the category of differential
 graded algebras. Clearly, one recovers the Hochschild cohomology of a
graded algebra $A$ by considering it as a differential graded algebra
with trivial differential. When $B=A$ with the canonical bimodule
structure, we write $\HHoch_{\mathsf {DG} }^\bullet(A)$ for 
$\HHoch_{\mathsf {DG} }^\bullet(A,A)$. As in the graded case, the differential
graded Hochschild complex $\Hoch_{\mathsf {DG} }(A)$ has a graded Lie algebra
structure, and both $\delta_0$ and $\delta_1$ are operators of adjoint
type for this Lie algebra structure.

\par
Since the vector field $Q$ squares to zero, it induces a differential
on the algebra of multivector fields of the differential graded
manifold $(N,Q)$. Algebraically, this amounts to saying that the
operator $\{\dd,\cdot\}$ acts as a differential on
$\calV^\bullet(A)$. We can therefore look at
$\calV^\bullet(A)$ as a bicomplex: the horizontal
differential is zero, and the vertical differential is
$\{\dd,\cdot\}$. We have a HKR map $\calV^\bullet(A)\to \Hoch_{\mathsf {DG} }(A)$, 
which is defined as in the case of differential algebras. 
\begin{Lem}
\label{lem:dr}
The HKR map $(\calV^\bullet(A),\{\dd,\cdot\},0)\to
(\Hoch_{\mathsf {DG} }(A),\delta_0,\delta_1)$ is a map of bicomplexes.
\end{Lem}
\begin{proof}
What we have said on the HKR map for graded algebras implies that
$\HKR\colon (\calV^\bullet(A),0)\to (\Hoch_{\mathsf {DG}
}(A),\delta_1)$ is a map of complexes. So we are left with checking the
compatibility of $\{\dd,\cdot\}$ with the differential
$\delta_0$. This follows from the following more general fact: given a
vector field $X$ and a multivector field $Y$, then
$\HKR(\{X,Y\})=\{\HKR(X),\HKR(Y)\}$, as one can easily verify. Note that
for an arbitrary multivector field $X$, the above identity only holds
up to homotopy. Since
$\delta_0(\HKR(Y))=\{\HKR(\dd),\HKR(Y)\}$, this concludes the proof.
\end{proof}

Being compatible with the differentials, the HKR map induces a map
between the cohomologies of the total complexes
$
{\mathsf H}^\bullet(\calV^\bullet(A),\{\dd,\cdot\})\to
\HHoch_{\mathsf {DG} }^\bullet(A),
$
which is a map of graded Lie algebras. In contrast with the case of
smooth algebras which are the subject of the classical HKR theorem, this
map is not an isomorphism in general, as the next theorem shows.

\begin{Thm}\label{thm:derham}
If $M$ is a simply connected closed oriented smooth manifold of
positive dimension, then the HKR map ${\mathsf H}^\bullet(\calV^\bullet
(\Omega^\bullet(M)),\{\dd,\cdot\})\to\HHoch_{\mathsf {DG} }^\bullet 
(\Omega^\bullet(M))$ is not an isomorphism.
\end{Thm}
We need the following Lemma, relating the $\{\dd,\cdot\}$-cohomology of
multivector fields on $T[1]M$ to the de Rham cohomology of $M$:
\begin{Lem}\label{lemma:mv-de-rham}
For any differential manifold $M$, there is an isomorphism 
\[
{\mathsf H}^\bullet(\calV^\bullet
(\Omega^\bullet(M),\{\dd,\cdot\}))\simeq {\mathsf H}^\bullet_{{\rm de
Rham}}(M).\]
\end{Lem}
\begin{proof}
Recall that $\calV^\bullet
(\Omega^\bullet(M))$ is the algebra of multivector fields on the graded
manifold $T[1]M$. We fix local 
coordinates $\{x^i,\theta^j\}$ on $T[1]M$, where $x^i$ are (even) 
coordinates on $M$ and $\theta^j$ (odd) coordinates on the fibers. 
Consider the globally well-defined
derivation $\iota_E$ which on the local generators of multivector fields acts
as
\begin{equation*}
\iota_E\left(x^i\right)=0\,;\quad
\iota_E\left(\theta^i\right)=0\,;\quad
\iota_E\left(\frac{\partial}{\partial x^i}\right)=\frac{\partial}{\partial 
\theta^i}\,;\quad
\iota_E\left(\frac{\partial}{\partial \theta^i}\right)=0\,.
\end{equation*}
The derivation $\{\dd,\cdot\}$ acts as
\begin{equation*}
\left\{\dd, x^i\right\}=\theta^i\,;\quad
\left\{\dd, \theta^i\right\}=0\,;\quad
\left\{\dd, \frac{\partial}{\partial x^i}\right\}=0\,;\quad
\left\{\dd, \frac{\partial}{\partial \theta^i}\right\}=\frac{\partial}
{\partial x^i}\,.
\end{equation*}
It follows that $L_E=\{\dd,\cdot\}\circ\iota_E + \iota_E\circ \{\dd,\cdot\}$ 
is a derivation on $\calV(T[1]M)$ which, when restricted
to the fields of degree $m$, is the multiplication by $m$; namely
\begin{equation*}
L_E\left(x^i\right)=0\,;\quad
L_E\left(\theta^i\right)=0\,;\quad
L_E\left(\frac{\partial}{\partial x^i}\right)=\frac{\partial}{\partial
x^i}\,;\quad L_E\left(\frac{\partial}{\partial
\theta^i}\right)=\frac{\partial}{\partial 
\theta^i}\,.
\end{equation*}
Now, suppose that $\Psi$ is a $\{\dd,\cdot\}$-closed multivector field of 
degree $m\ge 1$. Then it is also $\{\dd,\cdot\}$-exact:
\begin{align*}
\Psi&=\frac{1}{m}L_E(\Psi)=\frac{1}{m}\{\dd,\iota_E\Psi\} +
\frac{1}{m}\iota_E(\{\dd,\Psi\})\\
&=\{\dd,\frac{1}{m}\iota_E\Psi\}
\end{align*}
This shows that higher cohomology groups vanish, and we are left to prove that
${\mathsf H}^0(\calV^\bullet(T[1]M),\{\dd,\cdot\})={\mathsf H}^\bullet_{{\rm de
Rham}}(M)$. To see this, just notice that the
$0$-vector fields on $T[1]M$ are the differential forms on $M$ and the action
of
$\{\dd,\cdot\}$ on $\calV^0(T[1]M)$ is precisely the action of  the de Rham
differential on
$\Omega^\bullet(M)$.
\end{proof}

\begin{proof}[Proof of Theorem \ref{thm:derham}]
Let $\calL M$ be the free loop space on $M$. On the one hand we have
Chen's isomorphism
\cite{Chen, GJP}
\[
{\mathsf H}_\bullet(\calL M)[\dim M]
\simeq
\HHoch_{\mathsf {DG} }^\bullet(\Omega^\bullet(M)).
\]
On the other hand, we have the isomorphism 
\[
{\mathsf H}^\bullet(\calV^\bullet
(\Omega^\bullet(M),\{\dd,\cdot\}))\simeq {\mathsf H}^\bullet_{{\rm de
Rham}}(M)
\]
from Lemma \ref{lemma:mv-de-rham}.
Finally, $
{\mathsf H}_\bullet(\calL M)[\dim M]
\not\simeq {\mathsf H}^\bullet_{{\rm de Rham}}(M)$
for any simply connected closed oriented smooth manifold of
positive dimension \cite{SV}.
\end{proof}

\begin{Rem}
Observe that
another way of proving 
Lemma \ref{lemma:mv-de-rham}
goes through the Gerstenhaber
isomorphism described in Remark \ref{rem:Fourier}.
In fact, it is not difficult to see that the image of the multivector field
$\dd$ under this isomorphism is the restriction to $A=\Gamma(S^\bullet TM)$ 
of the canonical Poisson bivector field on the symplectic manifold $T^*M$.
Thus, ${\mathsf H}^\bullet(\calV^\bullet(T[1]M),\{\dd,\cdot\})$ is
isomorphic to the Poisson cohomology of $T^*M$ (restricted to functions
polynomial along the fibers) 
which in turn (by nondegeneracy of the Poisson structure) is isomorphic to 
the de Rham cohomology of the total space
and hence of the base.
\end{Rem}

\section{BV structures in the differential graded case}\label{sec:sq}

By forgetting the differential, i.e., by looking at a differential
graded manifolds simply as a graded manifolds, we obtain a BV
structure on the space of their multivector fields, as in Section
\ref{sec:mvf}. In general, this BV structure does not induce a BV
structure on the
$\{Q,\cdot\}$-cohomology of multivector fields. Indeed, the BV generator 
$\Delta$ is a derivation of the BV bracket, so it does not map 
$\{Q,\cdot\}$-closed vector fields to $\{Q,\cdot\}$-closed vector fields. 
Rather, if $X$ is a $\{Q,\cdot\}$-closed vector field, then
\[
\{Q,\Delta(X)\}=\{\Delta(Q),X\}.
\]
Yet, this implies that, if the vector field $Q$ is divergence-free,
i.e., if $\Delta(Q)=0$ then $\Delta$ induces a BV structure on the
$\{Q,\cdot\}$-cohomology, since
\[
\Delta\{Q,X\}=-\{Q,\Delta(X)\}
\]
and so $\{Q,\cdot\}$-exact multivector fields are mapped to 
$\{Q,\cdot\}$-exact multivector
fields. Note that, since the divergence operator $\Delta$ we are
considering in this paper is defined as the variation of the
Berezinian volume form of $N$ along a vector field, the condition
$\Delta(Q)=0$ means that the volume form is $Q$-invariant. 
A differential graded manifold $(N,Q)$ with a $Q$-invariant Berezinian
volume form is usually called an SQ-manifold \cite{schwarz,sch2}.
\begin{Rem}
In case $N$ is an odd symplectic manifold and the vector field $Q$ is 
Hamiltonian, one speaks of PQ-manifolds \cite{AKSZ}. Note that, 
if $Q=H_S$, i.e., if $S$ is the
function on $N$ whose Hamiltonian vector field is $Q$, then $\div(Q)=\Delta(S)$
and $\{Q,Q\}=H_{\{S,S\}}$ where on the right we have the odd Poisson bracket 
associated to the odd symplectic structure on $N$. Therefore, under the mild
assumption that $S$ has at least one critical point, the two equations 
$\{Q,Q\}=0$ and $\div(Q)=0$ 
imply the quantum master equation for
$S$, namely $\Delta(S)+\frac{\sqrt{-1}}{2\hbar}\{S,S\}=0.$ 
\end{Rem}
As far as concerns the BV structures on Hochschild cohomology, the 
same construction we described in Section~\ref{sec:hoch} also works 
in the differential graded case: if $(A,\dd)$ is the differential 
graded algebra of functions on the SQ-manifold $(N,Q)$, then we have 
a BV structure on $\HHoch_{\mathsf {DG} }^\bullet(A)$ under the 
hypothesis that $\wp\colon \Hoch_{\mathsf {DG} }^\bullet(A) \to 
\Hoch_{\mathsf {DG} }^\bullet(A,A^*)$ is a quasi-isomorphism. 
Moreover, by the same argument used in Section
\ref{sec:hkr}, the HKR map ${\mathsf
H}^\bullet(\calV^\bullet(A),\{\dd,\cdot\})
\to \Hoch^\bullet_{\mathsf {DG} }(A)$ is a BV map in this case.
\par
An example is given by the de Rham algebra $(\Omega^\bullet(M),\dd)$ of
a smooth closed manifold $M$. In the coordinates $\{x^i,\theta^j\}$ on
$T[1]M$, the de Rham differential on $\Omega^\bullet(M)$ is written
\[
\dd=\sum_{i=1}^{\dim M}\theta^i\frac{\partial}{\partial x^i},
\]
so that its divergence is
\[
\div(\dd)=\sum_{i=1}^{\dim M}\frac{\partial \theta^i}{\partial
x^i}=0.
\]
The pairing on $\Omega^\bullet(M)$ induced by the canonical Berezinian
volume form on $T[1]M$ is the usual Poincar\'e duality pairing:
\[
\langle\omega_1 ,\omega_2 \rangle=\int_M\omega_1\wedge\omega_2.
\] 
The induced map $\wp\colon 
\Hoch_{\mathsf {DG} }^\bullet(\Omega^\bullet(M) ) \to 
\Hoch_{\mathsf {DG} }^\bullet(\Omega^\bullet(M), \Omega^\bullet(M)^*)$ is a 
quasi-isomor\-phism \cite{M}, and so there exists 
a BV algebra structure on $\HHoch_{\mathsf {DG} }^\bullet
(\Omega^\bullet(M))$. This BV algebra structure coincides,
via Chen's isomorphism
\[
\HHoch_{\mathsf {DG} }^\bullet(\Omega^\bullet(M))\simeq
{\mathsf H}_\bullet(\calL M)[\dim M],
\]
 with the Chas--Sullivan BV structure on the
homology of the free loop space of
$M$ \cite{CS, Chen, CJ, GJP, M, Tr}. Also the $\{\dd,\cdot\}$-cohomology
of
$\calV^\bullet(\Omega^\bullet(M))$  has a nice geometrical
interpretation: we have shown in the proof of Lemma
\ref{lemma:mv-de-rham} that 
\[
{\mathsf H}^{p}(\calV^\bullet(\Omega^\bullet(M)),\{\dd,\cdot\}) =
\begin{cases}
0 & \text{ if }p\neq 0\\
{\mathsf H}^\bullet_{{\rm de
Rham}}(M) & \text{ if }p=0.
\end{cases}
\]
Note that, since the $\{\dd,\cdot\}$-cohomology
of
$\calV^\bullet(\Omega^\bullet(M))$ is concentrated in degree zero, the
BV structure on ${\mathsf
H}^\bullet(\calV^\bullet(\Omega^\bullet(M)),\{\dd,\cdot\})$ is the
trivial one. Finally, the BV map ${\mathsf
H}^\bullet(\calV^\bullet(\Omega^\bullet(M)),\{\dd,\cdot\})\to
\HHoch_{\mathsf {DG} }^\bullet(\Omega^\bullet(M))$ is the natural map
\[
{\mathsf H}^\bullet_{{\rm de
Rham}}(M)\simeq {\mathsf H}_\bullet(M)[\dim M]\to{\mathsf H}_\bullet(\calL M)[\dim M]
\]
induced by the natural embedding $M\hookrightarrow {\calL }M$ which
identifies the points of $M$ with the constant loops in ${\calL} M$.

\begin{Rem}
The constructions of Section \ref{sec:mdo} also work in the differential
graded case: a BV structure is defined on the total cohomology of any
sub-Gerstenhaber algebra $C_{\mathsf {DG} }^\bullet(A)$ of $\Hoch_{\mathsf {DG} }^\bullet(A)$, whose $\wp$-image in $\Hoch_{\mathsf {DG} }^\bullet(A,A^*)$ is
closed under $\beta$. This way we obtain a BV structure on the total
cohomology of multidifferential operators on an SQ-manifold. Moreover,
the HKR map ${\mathsf H}^\bullet(\calV^\bullet(A),\{\dd,\cdot\})\to {\mathsf H}_{\mathsf {DG} }^\bullet (\calD^\bullet(A))$ is a BV map.
\end{Rem}

\end{document}